\documentclass[12pt, eqno]{article}

\usepackage{amssymb,latexsym}
\usepackage{amsmath,latexsym}
\usepackage[top=1in,bottom=1in,left=1.3in,right=1.3in]{geometry}
\usepackage{amscd}
\usepackage{amsthm}

\newtheorem{theorem}{Theorem}[section]

\newtheorem{corollary}[theorem]{Corollary}

\newtheorem{example}[theorem]{Example}
\newtheorem{remark}[theorem]{Remark}
\newtheorem{lemma}[theorem]{Lemma}

\newtheorem{proposition}[theorem]{Proposition}

\begin{document}

\title{Notes on Ricci solitons in $f$-cosymplectic manifolds}

\author{\\  Xiaomin Chen
\thanks{
The author is supported by the Science Foundation of China
University of Petroleum-Beijing(No.2462015YQ0604) and partially supported
by  the Personnel Training and Academic
Development Fund (2462015QZDX02).
 }\\
{\normalsize College of  Sciences, China University of Petroleum (Beijing),}\\
{\normalsize Beijing, 102249, China}\\
{\normalsize xmchen@cup.edu.cn}}
\maketitle \vspace{-0.1in}


\abstract{The purpose of this article is to study an $f$-cosymplectic manifold $M$ admitting Ricci solitons.
Here we consider mainly two classes of Ricci solitons on $f$-cosymplectic manifolds. One is the class of contact Ricci solitons.
 The other is the class of gradient Ricci solitons,
for which we give the local classifications of
$M$. Meanwhile, we also give some properties of $f$-cosymplectic manifolds.}
 \vspace{-0.1in}
\medskip\vspace{12mm}

  \vspace{2mm}

 \vspace{2mm}

\section{Introduction}

In contact geometry, one important class of almost contact manifolds are almost Kenmotsu manifolds,
 which were introduced firstly by Kenmotsu in \cite{K}. Given an almost Kenmotsu structure $(\phi,\xi,\eta,g)$, we can
get an almost $\alpha$-Kenmotsu structure by a homothetic deformation:
\begin{equation*}
 \phi'=\phi,\eta'=\frac{1}{\alpha}\eta,\xi'=\alpha\xi,g'=\frac{1}{\alpha^2}
\end{equation*}
 for some non-zero real constant $\alpha$. Note that almost $\alpha$-Kenmotsu structures are related to
some special conformal deformations of almost cosymplectic structures (\cite{V}).

The notion of almost cosymplectic manifolds was first given by Goldberg and Yano in \cite{GY}. Later Kim and Pak in \cite{KP} defined
a new class called as {almost $\alpha$-cosymplectic manifolds} by combining almost cosymplectic and almost $\alpha$-Kenmotsu manifolds,
where $\alpha$ is a real number. Recently, Based on Kim and Pak's work, Aktan et al.\cite{AYM} considered a wide subclass of almost
contact manifolds, which are called {\it almost $f$-cosymplectic manifolds} defined by choosing a smooth function $f$ in the conception of almost $\alpha$-cosymplectic manifolds instead of any real number $\alpha$.

In the following we recall that a \emph{Ricci soliton} $(g,V)$ is a Riemannian metric $g$ together a vector field $V$ that satifies
\begin{equation}\label{1}
 \frac{1}{2}\mathcal{L}_V g+Ric-\lambda g=0,
\end{equation}
where $\lambda$ is constant and $V$ is called potential vector field.  The Ricci soliton is said to be {\it shrinking, steady and expanding} according
as $\lambda$ is positive, zero and negative respectively.  Specially, when the potential vector $V$ is taken as
the Reeb vector field on an almost contact metric manifold, it is called {\it contact Ricci soliton}, and if $V=DF$, the gradient vector
field of some function $F$ on $M$,  the Ricci soliton is  called accordingly a {\it gradient Ricci soliton}.
 The Ricci soliton is important not only for studying topology of
  manifold but also in the string theory. Compact Ricci solitons are the
fixed point of the Ricci flow: $\frac{\partial}{\partial t}g=-2Ric$
projected from the space of metrics onto its quotient modulo diffeomorphisms and scalings, and often
arise as blow-up limits for the Ricci flows on compact manifolds.

Concerning the study of the Ricci solitons
it has a long history and a lot of conclusions, were acquired, see \cite{Cao,CK,GS,GSC,MS,P,PW}etc.
 In particular, we note that Ghosh \cite{G} studied a three-dimensional Kenmotsu manifold admitting a Ricci soliton and proved it is of constant curvature $-1$.

As the generalization of Kenmotsu manifolds, in this paper, we study a normal almost $f$-cosymplectic manifold, which will be said to be an {\it $f$-cosymplectic manifold}, and get the classifications of $f$-cosymplectic manifolds whose metrics are contact Ricci solitons and gradient Ricci solitons, respectively. In order to prove our theorems we need some basic conceptions, which are presented in Section 2, and the main results and proofs are given in Section 3.

\section{Some basic conceptions and related results}

 Let $M^{2n+1}$ be a $(2n+1)$-dimensional Riemannian manifold.
An \emph{almost contact structure} on $M$ is a triple $(\phi,\xi,\eta)$, where $\phi$ is a
$(1,1)$-tensor field, $\xi$ a unit vector field, $\eta$ a one-form dual to $\xi$ satisfying
\begin{equation}\label{2.1}
\phi^2=-I+\eta\otimes\xi,\,\eta\circ \phi=0,\,\phi\circ\xi=0.
\end{equation}

A smooth manifold with such a structure is called an \emph{almost contact manifold}.
It is well-known that there exists a Riemannian metric $g$ such that
\begin{equation}\label{2.2}
g(\phi X,\phi Y)=g(X,Y)-\eta(X)\eta(Y),
\end{equation} for any $X,Y\in\mathfrak{X}(M)$. It is easy to get from \eqref{2.1} and \eqref{2.2} that
\begin{equation}\label{2.3}
g(\phi X,Y)=-g(X,\phi Y),\,g(X,\xi)=\eta(X).
\end{equation}
 An almost contact structure $(\phi,\xi,\eta)$ is said
to be \emph{normal} if the Nijenhuis torsion
\begin{equation*}
  N_\phi(X,Y)=\phi^2[X,Y]+[\phi X,\phi Y]-\phi[\phi X,Y]-\phi[X,\phi Y]+2d\eta(X,Y)\xi,
\end{equation*}
vanishes for any vector fields $X,Y$ on $M$.

Denote $\omega$ by the fundamental 2-form on $M$ defined by $\omega(X,Y):=g(\phi X,Y)$ for all $X,Y\in\mathfrak{X}(M)$.
If $\eta$ and $\omega$ are closed, then an almost contact structure is called {\it almost cosymplectic}, and it is said to be {\it cosymplectic} if
 in addition the almost contact structure is normal. An almost contact structure is said to be almost {\it $\alpha$-Kenmotsu} if $d\eta=0$ and $d\omega=2\alpha\eta\wedge\omega$ for a non-zero constant $\alpha$. More generally, if the constant $\alpha$ is any real number, then an almost contact
 structure is said to be
 {\it almost $\alpha$-cosymplectic}(\cite{OAM}). Moreover, Aktan et al.\cite{AYM} generalized the real number $\alpha$ to any smooth function $f$ on $M$ and defined an \emph{almost $f$-cosymplectic manifold}, which is an almost contact metric manifold $(M,\phi,\xi,\eta,g)$ such that
 $d\omega=2f\eta\wedge\omega$ and $d\eta=0$ for a smooth function $f$ satisfying $df\wedge\eta=0$.

In addition, if the almost $f$-cosymplectic structure on $M$ is
 normal, we say that $M$ is an \emph{$f$-cosymplectic manifold.} Obviously, if $f$ is constant, then an $f$-cosymplectic
manifold is either cosymplectic under condition $f=0$, or $\alpha$-Kenmotsu ($\alpha=f\neq0$). Furthermore, there exists a distribution
$\mathcal{D}$ of an $f$-cosymplectic manifold is defined by $\mathcal{D}=\ker\eta$,  which is integrable since $d\eta=0$.

Besides, for an almost contact manifold $(M,\phi,\xi,\eta,g)$,
we denote $h:=\frac{1}{2}\mathcal{L}_\xi\phi$, which is a self-dual operator.
Since an $f$-cosymplectic manifold is normal, $h=0$. Therefore, in virtue of \cite[Proposition 9, Proposition 10]{AYM} we know that
for a $(2n+1)$-dimensional $f$-cosymplectic manifold the following identities are valid:
 \begin{align}
\nabla_X\xi=&-f\phi^2X,\label{2.4}\\
  Q\xi =&-2n\widetilde{f}\xi,\label{2.5} \\
   R(X,Y)\xi=&\widetilde{f}[\eta(X)Y-\eta(Y)X],\label{2.6}
\end{align}
where $\nabla$ and $Q$ denote respectively the Levi-Civita connection and Ricci operator of $M$, and $\widetilde{f}\triangleq\xi(f)+f^2$.

\begin{proposition}\label{p2.1}
For any an  $f$-cosymplectic manifold, if $\xi(\widetilde{f})=0$, then $\widetilde{f}=const.$
\end{proposition}
\proof  Differentiating \eqref{2.6} along any vector field $Z$ we have
\begin{align*}
(\nabla_ZR)(X,Y)\xi=&\nabla_Z(R(X,Y)\xi)-R(\nabla_ZX,Y)\xi-R(X,\nabla_ZY)\xi-R(X,Y)\nabla_Z\xi\\
=&Z(\widetilde{f})[\eta(X)Y-\eta(Y)X]+\widetilde{f}f[g(X,Z)Y-g(Y,Z)X]\\
&-fR(X,Y)Z.
\end{align*}
Then using the second Bianchi identity
\begin{equation*}
 (\nabla_ZR)(X,Y)\xi+(\nabla_XR)(Y,Z)\xi+(\nabla_YR)(Z,X)\xi=0,
\end{equation*}
 we have
\begin{align*}
[&Y(\widetilde{f})\eta(Z)-Z(\widetilde{f})\eta(Y)]X+[Z(\widetilde{f})\eta(X)-X(\widetilde{f})\eta(Z)]Y\\
&+[X(\widetilde{f})\eta(Y)-Y(\widetilde{f})\eta(X)]Z-f[R(X,Y)Z+R(Y,Z)X+R(Z,X)Y]=0.
\end{align*}
By taking $Z=\xi$ and using \eqref{2.6}, we know
\begin{equation}\label{2.7}
\xi(\widetilde{f})[\eta(Y)X-\eta(X)Y]-X(\widetilde{f})\phi^2Y+Y(\widetilde{f})\phi^2X=0.
\end{equation}

If we assume $\xi(\widetilde{f})=0$,
then we can obtain $X(\widetilde{f})=0$ for every vector field $X$ by taking the inner product of \eqref{2.7} with $Y$, putting $Y=e_i$ and summing over $i$ in the resulting equation (where $\{e_i\}$ is the local orthonormal frame of $M$). \qed\bigskip

Obviously, it reduces directly the following corollary in view of $df\wedge\eta=0$.
\begin{corollary}
An $f$-cosymplectic manifold is a cosymplectic manifold if $f$ vanishes along $\xi$.
\end{corollary}
\begin{proposition}\label{p2.3*}
A compact $f$-cosymplectic manifold $M^{2n+1}$ with $\xi(\widetilde{f})=0$ is $\alpha$-cosymplectic. In particular, if $\widetilde{f}=0$, $M$ is cosymplecitc.
\end{proposition}
\proof  As $\xi(\widetilde{f})=\xi(\xi(f))+2f\xi(f)=0$, we obtain $\xi(\xi(f))=-2f\xi(f)$. On the other hand, we know that $f$ satisfies $df\wedge\eta=0$,
that means that $Df=\xi(f)\xi$, where $D$ is the gradient operator with respect to $g$. For every field $X$, it follows from \eqref{2.4} that
\begin{equation*}
  \nabla_XDf=X(\xi(f))\xi+\xi(f)\nabla_X\xi=X(\xi(f))\xi-f\xi(f)\phi^2X.
\end{equation*}

Since $\nabla_\xi\xi=0$, for every point $p\in M$, we may take a locally orthonormal basis $\{e_i\}$ of $T_pM$ such that $e_{2n+1}=\xi$ and $\nabla_{e_i}e_i=0$.
Therefore
\begin{equation}\label{2.8}
\Delta f=\sum_ig(\nabla_{e_i}Df,e_i)=\xi(\xi(f))+2nf\xi(f)=(1-n)\xi(\xi(f)),
\end{equation}
 where $\Delta$ is the Laplace operator. Since $e_i(f)=g(Df,e_i)=0$ for $i=1,\cdots,2n$, we find $\Delta f=\sum_ie_i(e_i(f))=\xi(\xi(f))$. Hence it yields from \eqref{2.8} that
$\xi(\xi(f))=0$, which shows $\Delta f=0$, i.e. $f$ is constant. If $\widetilde{f}=0$, i.e. $0=\xi(f)+f^2=f^2$, then $f=0$.\qed
\begin{remark}\label{r1}
In \cite{B},  Blair proved that a cosymplectic manifold is locally the product of a K\"{a}hler manifold and an interval or unit circle $S^1$.
\end{remark}

Moreover, for the case of three-dimension, we have
\begin{lemma}\label{L1}
For a three-dimensional $f$-cosymplectic manifold $M^3$, we have
\begin{equation}\label{2.9}
  QY=\Big(-3\widetilde{f}-\frac{R}{2}\Big)\eta(Y)\xi+\Big(\widetilde{f}+\frac{R}{2}\Big)Y,
\end{equation}
where $R$ is the scalar curvature of $M$.
\end{lemma}
\proof It is well-known that the curvature tensor of any three-dimensional Riemannian manifold is written as
\begin{align}\label{2.10}
  R(X,Y)Z =&g(Y,Z)QX-g(X,Z)QY+Ric(Y,Z)X \nonumber\\
         &-Ric(X,Z)Y-\frac{R}{2}[g(Y,Z)X-g(X,Z)Y].
\end{align}

Putting $Z=\xi$ and using \eqref{2.5}, \eqref{2.6}, we have
\begin{equation*}
  \Big(\widetilde{f}+\frac{R}{2}\Big)\Big(\eta(Y)X-\eta(X)Y\Big)=\eta(Y)QX-\eta(X)QY.
\end{equation*}
Moreover, by taking $X=\xi$ and using \eqref{2.6} again, we obtain \eqref{2.9}.\qed

\section{Main results and proofs}
In this section we mainly discuss two classes of Ricci solitons, i.e., contact Ricci solitons and gradient Ricci solitons
in $f$-cosymplectic manifolds, respectively. At first, for a general Ricci soliton we have the following lemma, which was showed
by Cho.
\begin{lemma}(\cite[Lemma 3.1]{Cho})\label{l3.1}
If $(g,V)$ is a Ricci soliton of a Riemannian manifold then we have
\begin{equation*}
  \frac{1}{2}||\mathcal{L}_Vg||^2=V(R)+2{\rm div}(\lambda V-QV),
\end{equation*}
where $R$ denotes the scalar curvature.
\end{lemma}
\begin{theorem}
If an $f$-cosymplectic manifold $M^{2n+1}$ admits a contact Ricci soliton, then $M^{2n+1}$
 is locally isometric to the product of a line and a Ricci-flat K\"{a}hler (Calabi-Yau) manifold.
\end{theorem}
\proof In view of \eqref{2.4}, we have
\begin{equation*}
  (\mathcal{L}_\xi g)(X,Y)=g(\nabla_X\xi,Y)+g(X,\nabla_Y\xi)=2f[g(X,Y)-\eta(X)\eta(Y)].
\end{equation*}
Therefore it implies from the Ricci equation \eqref{1} with $V=\xi$ that
\begin{equation}\label{3.1}
  Ric(X,Y)=(\lambda-f)g(X,Y)+f\eta(X)\eta(Y).
\end{equation}
By \eqref{3.1}, the Ricci operator $Q$ is provided
\begin{equation*}
 QX=(\lambda-f)X+f\eta(X)\xi.
\end{equation*}
for any vector field $X$ on $M$. Thus
\begin{align}
Q\xi=&\lambda\xi,\label{3.2}\\
R=&(2n+1)\lambda-2nf.\label{3.3}
\end{align}
By Lemma \ref{l3.1}, \eqref{3.2} and \eqref{3.3}, we find that
\begin{equation}\label{3.4}
  \frac{1}{2}||\mathcal{L}_\xi g||^2=-2n\xi(f).
\end{equation}

Since $(\mathcal{L}_\xi g)(X,Y)=2fg(\phi X,\phi Y)$ and $\widetilde{f}=-\frac{\lambda}{2n}$ is constant
  followed by comparing \eqref{2.5} with \eqref{3.2}, a straightforward computation implies
$f^2=const.,$ i.e., $f$ is constant. Hence $\xi$ is Killing by \eqref{3.4}.
Moreover we get $f=0$ from $(\mathcal{L}_\xi g)(X,Y)=2fg(\phi X,\phi Y)$.
Namely, $M$ is cosymplectic. Further we have $Ric=0$ since $\lambda=-2n\widetilde{f}=0$.
Therefore we complete the proof of our result. \qed\bigskip

In view of the above proof, we get immediately the following corollary.
\begin{corollary}
A contact Ricci soliton in an $f$-cosymplectic manifold is steady.
\end{corollary}

In the following we assume that  an $f$-cosymplectic manifold $M^{2n+1}$ admits a gradient Ricci soliton
and the function $f$ satisfies $\xi(\widetilde{f})=0$.
\begin{theorem}\label{T1}
Let $M^{2n+1}$ be an $f$-cosymplectic manifold with a gradient Ricci soliton.
If $\xi(\widetilde{f})=0$, then one of the following statements holds:
\begin{enumerate}
  \item $M$ is locally the product of a K\"{a}hler manifold and an interval or unit circle $S^1$,
  \item  $M$ is Einstein.
\end{enumerate}
\end{theorem}

In order to prove the theorem  we first prove
\begin{lemma}\label{p3.7}
Let $M^3$ be a three-dimensional $f$-cosymplectic manifold with a Ricci soliton.
Then the following equation holds:
\begin{equation}\label{3.5}
2\xi(\widetilde{f})+\frac{\xi(R)}{2}+2(3\widetilde{f}+\frac{R}{2})f=0.
\end{equation}
\end{lemma}
\proof By Lemma \ref{L1}, we compute
\begin{equation}\label{3.6}
\begin{aligned}
  (\nabla_XRic)(Y,Z)=&[-3X(\widetilde{f})-\frac{X(R)}{2}]\eta(Y)\eta(Z)+[-3\widetilde{f}-\frac{R}{2}](\nabla_X\eta)(Y)\eta(Z)\\
  &+[-3\widetilde{f}-\frac{R}{2}]\eta(Y)(\nabla_X\eta)(Z)+[X(\widetilde{f})+\frac{X(R)}{2}]g(Y,Z)\\
  =&[-3X(\widetilde{f})-\frac{X(R)}{2}]\eta(Y)\eta(Z)+[-3\widetilde{f}-\frac{R}{2}]fg(\phi X,\phi Y)\eta(Z)\\
  &+[-3\widetilde{f}-\frac{R}{2}]f\eta(Y)g(\phi X,\phi Z)+[X(\widetilde{f})+\frac{X(R)}{2}]g(Y,Z).
\end{aligned}
\end{equation}
 Notice that for every vector $Z$,
the following relation holds:
\begin{align}\label{3.7}
  \sum_{i=1}^3\Big[(\nabla_ZRic)(e_i,e_i)-2(\nabla_{e_i}Ric)(e_i,Z)\Big]=0,
  \end{align}
which is followed from the formulas (8) and (9) of \cite{G}, where $\{e_1,,e_{2},e_{3}=\xi\}$ is a local orthonormal frame of $M$.

Making use of \eqref{3.6}, we obtain from \eqref{3.7} that
\begin{align*}
\Big[-3\xi(\widetilde{f})-\frac{\xi(R)}{2}\Big]\eta(Z)+2(-3\widetilde{f}-\frac{R}{2})f\eta(Z)+Z(\widetilde{f})=0.
  \end{align*}
Putting $Z=\xi$ in the above formula gives \eqref{3.5}. \qed\bigskip\\
{\it Proof of Theorem \ref{T1}} By Proposition \ref{p2.1}, $\widetilde{f}=constant.$ It is clear that the
Ricci soliton equation \eqref{1} with $V=DF$ for some smooth function $F$ implies
\begin{equation}\label{3.8}
  \nabla_YDF=-QY+\lambda Y.
\end{equation}
Therefore we have $R(X,Y)DF=(\nabla_YQ)X-(\nabla_XQ)Y$. Putting $Y=\xi$ further gives
\begin{equation}\label{3.9}
  R(X,\xi)DF=(\nabla_\xi Q)X-(\nabla_XQ)\xi.
\end{equation}

On the other hand, from \eqref{2.6} and the Bianchi identity, we have
\begin{equation}\label{3.10}
  R(X,\xi)Y=\widetilde{f}[g(X,Y)\xi-\eta(Y)X].
\end{equation}
Replacing $Y$ by $DF$ in \eqref{3.10} and comparing with \eqref{3.9}, we get
\begin{equation}\label{3.11}
  (\nabla_\xi Q)X-(\nabla_XQ)\xi=\widetilde{f}[X(F)\xi-\xi(F)X].
\end{equation}
Taking the inner product of the previous equation with $\xi$ and using \eqref{2.5}, we arrive at that
\begin{equation}\label{3.12}
\widetilde{f}[X(F)-\xi(F)\eta(X)]=0.
\end{equation}

Next we divide into the following cases.

{\bf Case I:} $\widetilde{f}=0$ and $n>1$. That is, $\xi(f)=-f^2$, then $Df=-f^2\xi.$ If $f\not\equiv0$
then there is
an open neighborhood $\mathcal{U}$ such $f|_\mathcal{U}\neq0$, thus in this case
$\xi=-\frac{Df}{f^2}=D(\frac{1}{f}).$ Since $\Delta f=0$(see the proof Proposition \ref{p2.3*}),
\begin{equation*}
  0=\Delta(f\cdot\frac{1}{f})=\frac{1}{f}\Delta f+2g(Df,D(\frac{1}{f}))+f\Delta\frac{1}{f}=2\xi(f)+f\text{div}\xi.
\end{equation*}
From \eqref{2.4}, we know $\text{div}\xi=2nf$. When $n>1$, substituting this
into the previous equation implies $f=0,$ which leads to a contradiction. Hence $f\equiv0$, that is, $M$ is cosymplectic.

{\bf Case II:}
 $\widetilde{f}\neq0$.  By \eqref{3.12}, the following identity is obvious:
\begin{equation}\label{3.13}
  DF=\xi(F)\xi.
\end{equation}
 Substituting this into \eqref{3.8} and using \eqref{2.4}, we give
\begin{equation}\label{3.14}
Y(\xi(F))\xi-f\xi(F)\phi^2Y=-QY+\lambda Y.
\end{equation}
 By taking an inner product with $\xi$ and using \eqref{2.5}, we further find
\begin{equation}\label{3.15}
 Y(\xi(F))=(2n\widetilde{f}+\lambda)\eta(Y).
\end{equation}
 Now taking \eqref{3.15} into \eqref{3.14} implies that for every vector $X$,
\begin{equation}\label{3.16}
  \lambda g(X,Y)-Ric(X,Y)=(2n\widetilde{f}+\lambda)\eta(X)\eta(Y)+f\xi(F)g(\phi X,\phi Y).
\end{equation}
Moreover, we derive from \eqref{3.16} that the scalar curvature
\begin{equation}\label{3.17}
  R=2n(-\widetilde{f}+\lambda-f\xi(F)).
\end{equation}

On the other hand, using \eqref{3.13} and \eqref{2.5}, we have
\begin{equation}\label{3.18}
  Ric(X,DF)=\xi(F)g(QX,\xi)=-2n\widetilde{f}\eta(X)\xi(F).
\end{equation}
It is well known that for any vector field $X$ on $M$,
\begin{equation}\label{3.19}
g(DR,X)=2Ric(DF,X),
\end{equation}
 which can be found in \cite{H}.
Applying \eqref{3.18} and \eqref{3.17} in this identity, we have
\begin{equation}\label{3.20}
X(f)\xi(F)+(2n\widetilde{f}+\lambda)\eta(X)=2f\xi(F)\eta(X).
\end{equation}
Putting $X=\xi$ the equation \eqref{3.20}, we get
\begin{equation}\label{3.21}
  (\xi(f)-2\widetilde{f})\xi(F)+f(2n\widetilde{f}+\lambda)=0.
\end{equation}
Differentiating \eqref{3.21} along $\xi$, we obtain from \eqref{3.15} that
\begin{equation}\label{3.22}
\xi(\xi(f))\xi(F)+2(\xi(f)-\widetilde{f})(2n\widetilde{f}+\lambda)=0.
\end{equation}
Since $\xi(\widetilde{f})=0$, we have $\xi(\xi(f))=-2f\xi(f)$. Substituting this into \eqref{3.22} yields
\begin{equation}\label{3.23}
  f\xi(f)\xi(F)+f^2(2n\widetilde{f}+\lambda)=0.
\end{equation}
Differentiating the above formula again along $\xi$, we obtain
\begin{equation*}
  (\xi(f)^2-2f^2\xi(f))\xi(F)+3f\xi(f)(2n\widetilde{f}+\lambda)=0.
\end{equation*}
Applying \eqref{3.21} in this equation implies
\begin{equation*}
  (\xi(f)+4f^2)\xi(f)\xi(F)=0.
\end{equation*}

Now if $\xi(f)+4f^2=0$ on some neighborhood $\mathcal{O}$ of $p\in M$, then $3f^2=-\widetilde{f}$ is constant, i.e., $f$ is constant on $\mathcal{O}$.
 Further we know $f=0$, which implies $\widetilde{f}=0$ on $\mathcal{O}$. It is a contradiction with the assumption $\widetilde{f}\neq0$.
 Therefore $\xi(f)\xi(F)=0$, and it follows from \eqref{3.23} that
 \begin{equation*}
   f^2(2n\widetilde{f}+\lambda)=0.
 \end{equation*}

If $2n\widetilde{f}+\lambda=0$, then it reduces from \eqref{3.21} that $(\xi(f)-2\widetilde{f})\xi(F)=0$, i.e., $(\xi(f)+2f^2)\xi(F)=0$.
As before we know $\xi(F)=0.$ It shows that $DF$ is identically zero because of \eqref{3.13}. Thus $M$ is Einstein.
Moreover, from \eqref{3.17} we get $R=2n(\lambda-\widetilde{f})$.

If $2n\widetilde{f}+\lambda\neq0$, we have $f\equiv0$, that is, $M$ is cosymplectic.

In particular, when $n=1$, we know that $\lambda+2\widetilde{f}=0$ and $R=2(\lambda-\widetilde{f})$. Moreover, using Lemma \ref{L1} we obtain
\begin{equation*}
QY=(-2\widetilde{f}-\lambda)\eta(Y)\xi+\lambda Y=-2\widetilde{f}Y.
\end{equation*}

{\bf Case III:} $\widetilde{f}=0$ and $n=1$. If $f\equiv0$, $M$ is cosymplectic. Next we always assume $f\neq0$ on some neighborhood.  By \eqref{2.5} we have $Q\xi=0$ when $\widetilde{f}=0$. Because $\widetilde{f}$ is constant, so we obtain $\xi(R)=0$ from \eqref{3.19}. That means $R=0$ by \eqref{3.5}.
 Moreover, in view of Lemma \ref{L1} we get $Q=0$.

Summarizing the above discussion,  we have proved that either $f\equiv0$ or $QY=-2\widetilde{f}Y$, thus by Remark \ref{r1} we complete the proof of theorem.\qed\\

Since an $\alpha$-cosymplectic manifold is actual an $f$-cosymplectic  manifold such that $f$ is constant, we obtain from Theorem \ref{T1}.
\begin{corollary}\label{C1}
Let $M^{2n+1}$ be an $\alpha$-cosymplectic manifold with a gradient Ricci soliton. Then $M$ is either locally the product of a K\"{a}hler manifold and an interval or unit circle $S^1$, or Einstein.
\end{corollary}

We note that Perelman in \cite{P} proved that on a compact Riemannian manifold a Ricci soliton is always a gradient Ricci soliton,
thus the following corollary is clear by Theorem \ref{T1}.
\begin{corollary}
Let $M^{2n+1}$ be a compact $f$-cosymplectic manifold with a Ricci soliton.
If $\xi(\widetilde{f})=0$, then $M$ is either locally the product of a K\"{a}hler manifold and an interval or unit circle $S^1$, or Einstein.
\end{corollary}

When $n=1$, we have
\begin{corollary}
Let $M^3$ be a three-dimensional $\alpha$-cosymplectic manifold with a Ricci soliton. If $\xi(R)=0$ then
 $M$ is either locally the product of a K\"{a}hler manifold and an interval or unit circle $S^1$, or Einstein.
 \end{corollary}
\proof Since an $\alpha$-cosymplectic manifold is an $f$-cosymplectic manifold with $f=\alpha$ is constant, we have $\widetilde{f}=\alpha^2$.
Since $\xi(R)=0$, making using of \eqref{3.5} we obtain
$$(3\alpha^2+\frac{R}{2})\alpha=0.$$
Therefore $\alpha=0$ or $R=-6\alpha^2$. We complete the proof by Lemma \ref{L1}.\qed\bigskip

Finally we give an example of an $f$-cosymplectic manifold satisfying $\xi(\widetilde{f})=0$.
\begin{example}
As the Example of three dimension in \cite{AYM} we also consider a three-dimensional manifold $M=\{(x,y,z)\in\mathbb{R}^3\}$, where $x,y,z$ are the standard coordinates in $\mathbb{R}^3$.
On $M$  we define the Riemannian metric
\begin{equation*}
  g=\frac{1}{e^{2\theta(z)}}(dx\otimes dx+dy\otimes dy)+dz\otimes dz,
\end{equation*} where $\theta(z)$ is a smooth function on $M$.
\end{example}

Clearly, the vector fields
\begin{equation*}
  e_1=e^{\theta(z)}\frac{\partial}{\partial x},\quad e_2=e^{\theta(z)}\frac{\partial}{\partial y},\quad e_3=\frac{\partial}{\partial z}
\end{equation*}
are linearly independent with respect to $g$ at each point of $M$. Also,  we see $g(e_i,e_j)=\delta_{ij}$ for $i,j=1,2,3$.

Let $\eta$ be the 1-form defined by $\eta(X)=g(X,e_3)$ for every field $X$  and $\phi$ be the $(1,1)$-tensor field defined by
$\phi(e_1)=e_2,\phi(e_2)=-e_1,\phi(e_3)=0$. Hence it is easy to get that $\eta=dz$, $\omega(e_1,e_2)=g(\phi(e_1),e_2)=1$ and $\omega(e_1,e_3)=\omega(e_2,e_3)=0$.

Furthermore, a straightforward computation gives the brackets of the vector fields $e_1,e_2,e_3$:
\begin{equation*}
  [e_1,e_2]=0,\quad [e_1,e_3]=-\theta'(z)e_1,\quad [e_2,e_3]=-\theta'(z)e_2.
\end{equation*}
Consequently, the the Nijenhuis torsion of $\phi$ is zero, i.e., $M$ is normal.

On the other hand, as in \cite{AYM}, it easily follows
\begin{equation*}
  \omega=\frac{1}{e^{2\theta(z)}}dx\wedge dy
\end{equation*}
and
\begin{equation*}
  d\omega=-2\theta'(z)e^{-2\theta(z)}dx\wedge dy\wedge dz=2\theta'(z)\omega\wedge\eta.
\end{equation*}
Therefore $M$ is an $f$-cosymplectic manifold with $f(x,y,z)=\theta'(z)$.

In order that $\xi(\widetilde{f})=0$, i.e.
$$
 \theta'''(z)+2\theta'(z)\theta''(z)=0,
$$
we need $ \theta''(z)+[\theta'(z)]^2=c$ for a constant $c$. In view of theory of ODE, the above equation is solvable.

\end{document}